\input moshe.fnt
\input moshe.mac
\input danny.mac
\input xy
\xyoption {all}

\null
\pageno=0
\centerline{ \footnote{\strut}{The second author was partially supported by the Israel Science Foundation (grant No. 343/07), and by an ERC grant.}{\titlefont
\uppercase{Admissible groups over two dimensional}}}
\centerline{{\titlefont
\uppercase{complete local domains}}}

\bigskip
\centerline{by}
\medskip
\centerline{\authorfont Danny Neftin}
\centerline{\affiliationfont
Department of Mathematics, Technion -- Israel Institute of Technology,}\centerline{\affiliationfont
Haifa 3200, Israel}

\centerline{\affiliationfont e-mail: neftind\at tx.technion.ac.il}

\medskip
\centerline{and}
\centerline{\authorfont Elad Paran}
\centerline{\affiliationfont
Einstein Institute of Mathematics, J.~Safra Campus, Givat Ram, }
\centerline{\affiliationfont
The Hebrew University of Jerusalem, Jerusalem 91904, Israel}
\centerline{\affiliationfont e-mail: paranela\at post.tau.ac.il}

\bigskip
\centerline{\smc Abstract}

{\narrower\smallskip\noindent
Let $K$ be the quotient field of
a complete local domain of dimension 2 with a separably closed residue field.
Let $G$ be a finite group of order not divisible by $\chr(K)$. Then $G$ is admissible over $K$ if and only if
its Sylow subgroups are abelian of rank at most $2$. \narrower\smallskip\noindent
\smallskip}

    \vskip1cm
\rightline{MR Classification: 12E30, 16S35}

\head Introduction \endhead


Complete local domains play an important role in commutative algebra and algebraic geometry, and their algebraic properties were already described in 1946 by Cohen's structure theorem. The Galois theoretic properties of their quotient fields were extensively studied over the past two decades. The pioneering work in this line of research is due to Harbater [Ha87], who introduced the method of patching to prove that if $R$ is a complete local domain with quotient field $K$, then every finite group occurs as a Galois group over $K(x)$. This result was strengthened by Pop [Po96], and in a different language, by Haran-Jarden [HaJ98], who showed that if $R$ is also of dimension 1, then every finite split embedding problem over $K(x)$ is solvable.

The first step towards higher dimension was made by Harbater-Stevenson [HaS05], who essentially showed that if $R$ is a complete local domain of dimension 2, then every finite split embedding problem over $\Quot(R)$ has $|R|$ independent solutions. That is, the absolute Galois group of $\Quot(R)$ is semi-free of rank $|R|$ (see [BHH10] for details on this notion). This result was later generalized by Pop [Po10] and by the second author [Pa10], who showed that if $K$ is the quotient field of a complete local domain of dimension exceeding 1, then $\Gal(K)$ is semi-free.
%

However, despite the major progress made in the study of Galois theory over these fields, little is known about the structure of division algebras over them. A step in that direction was recently made by Harbater-Hartmann-Krashen [HHK10]. In that work, the authors consider a question relating Galois theory and Brauer theory over a field $F$ -- which groups are admissible over $F$ ? That is, which finite groups occur as a Galois group of an {\bf adequate} Galois extension $L/F$ (recall that an extension $L/F$ is called adequate if $L$ is a maximal subfield in an $F$-central division algebra). Equivalently, for which groups $G$ there is a $G$-crossed product division algebra with center $F$.
Note that for $F$ as above and a finite extension $L/F$, the above maximality requirement can be omitted since any $L$ which is a subfield of an $F$-division algebra is also a maximal subfield of some $F$-division algebra (see Remark \cp2.9).

This question was first considered by Schacher [Sch68]. In that work, Schacher proved that any $\bbQ$-admissible group has metacyclic Sylow subgroups and conjectured the converse. Admissibility was studied extensively over global fields (for example: [Sch68], [Ste82], [Son83], [SS92], [Lid96], [Fei04], [Nef10] , function fields and fields of Laurant series (for example: [FSS92], [FS95a], [FS95b]).

The main theorem of [HHK10] asserts that if $F$ is a function field of one variable over a complete valued field with an algebraically closed residue field, then a finite group $G$, of order not divisible by $\chr(F)$, is admissible over $F$ if and only each of the Sylow subgroups of $G$ is abelian of rank at most $2$ (i.e. is generated by two elements).

In this work, we take the next natural step, and determine the admissible groups over quotient fields of two dimensional complete local domains, with a separably closed residue field. This problem was posed to the first author by Harbater. We show that the result of [HHK10] holds over these fields as well.

\proclaim{Main Theorem} Let $R$ be a complete local domain of dimension 2, with a separably closed residue field. Let $G$ be a finite group of order not divisible by $\chr(R)$. Then $G$ is admissible over $\Quot(R)$ if and only if each of its Sylow subgroups is abelian of rank at most 2. \endproclaim

The forward direction of our main theorem is essentially proven in [HHK10], using results of Colliot-Th\'el\`ene-Ojanguren-Parimala [COP02]. Our proof strategy for the backward direction is as follows. We first use Cohen's structure theorem to reduce the problem from $\Quot(R)$ to a field $E$ of the form $C((X,Y))$, where $C$ is a separably closed field. We then apply a patching argument as in [HHK10]; We explicitly realize each Sylow subgroup of $G$ by a Galois extension of $E$ which is a maximal subfield in some $E$-central division algebra. We then patch these realizations into a Galois extension $F/E$ with group $G$, in a way that also patches the division algebras into an $E$-central division algebra $D$ containing $F$ as maximal subfield.

The main innovation of this work, and the key ingredient of the proof, is a patching machinery over fields of the form  $K((X,Y))$. In [Po10] and [Pa10], problems over $K((X,Y))$ are lifted to $K((X,Y))(Z)$, solved there (via different methods), and then the solutions are specialized into solutions over $K((X,Y))$ using Hilbertianity. This approach seems inapplicable to our current goal, since adequate extensions usually do not remain adequate under specialization. Instead we patch groups directly over $K((X,Y))$. To this end we define "analytic fields" over $K((X,Y))$, satisfying the axioms of algebraic patching (i.e. matrix factorization and intersection), as presented in [HaJ98] (a slightly different axiomatization from the "field patching" axiomatization of [HaH09]). The construction of these analytic fields is based on a mixture of ideas from formal geometry and rigid analytic geometry, which we develop ad-hoc here.

\head \cs1. Analytic fields \endhead

In this section we establish our patching machinery.
Fix an infinite field $K$, and let $E = K((X,Y)) = \Quot(K[[X,Y]])$ be the field of formal power series over $K$ in $X,Y$. Denote by $v$ the order function of the maximal ideal $\langle X,Y \rangle$ in $K[[X,Y]]$, and extend $v$ (in the unique possible way) to a discrete valuation of $E$.

\reminfo{Construction \cp1.1}{Analytic rings over $E$} Let $I$ be a finite set. For each $i \in I$ let $c_i \in K$, such that $c_i \neq c_j$ for $i \neq j$ (such a choice is possible since $K$ is infinite). For each $i \in I$ denote $z_i = {Y \over X - c_iY} \in E$. For each $J \subseteq I$, consider the subring $K[z_j \st j \in J][X,Y]$ of $E$, and let $D_J$ be the completion of this ring with respect to $v$. Note that for each $J \subseteq I$, $D_J \subseteq D_I$, and that $D_\emptyset = K[[X,Y]]$, since $K[[X,Y]]$ is complete. Let $Q = \Quot(D_I)$, and for each $i \in I$ denote $Q_i = E \cdot D_{I \smallsetminus \{i\}}$, and $Q_i' = \bigcap_{j \neq i}Q_j$. \endrem

For the rest of this section, we fix the notation of Construction \cp1.1.

\proclaim{Lemma \cp1.2} Let $i \in I$. Then $v$ is trivial on $K(z_i)$. \endproclaim \demo{Proof} It suffices to prove that $v$ is trivial on $K[z_i]$. Let $0 \neq p = \sum_{n=0}^d a_n z_i^n \in K[z_i]$, with $a_0,\ldots,a_d \in K$. Without loss of generality, $a_0 \neq 0$. We have

$$\sum_{n=0}^d a_n z_i^n = {\sum_{n=0}^d a_n Y^n (X - c_iY)^{d-n} \over (X-c_iY)^d}.$$ By opening parentheses, the numerator in this expression can be written as a sum of monomials of degree $d$. For $n = 0$ we get the summand $a_0 X^d$, while all other monomials in this presentation have a positive power of $y$, so they do not cancel $a_0X^d$. Thus the numerator has value $d$, and clearly so does the denominator, so $v(p) = 0$. \qed

\proclaim{Corollary \cp1.3} The valuation $v$ is trivial on $K[z_i \st i \in I]$. \endproclaim \demo{Proof} Note that for each $i,j \in I$ we have $K(z_i) = K(z_j)$. Thus by the preceding lemma, $v$ is trivial on $K(z_i) = K(z_j \st j \in I)$, and in particular $v$ is trivial on the subring $K[z_i \st i \in I]$. \qed

\proclaim{Lemma \cp1.4} Let $J \subseteq I$ and $j \in J$. Then the ring $K[z_l \st l \in J][X - c_jY]$ is isomorphic to the ring of polynomials in one variable over $K[z_l \st l \in J]$. \endproclaim \demo{Proof} To prove the claim we show that if $\sum_{n = 0}^d a_n(X-c_jY)^n = 0$, for $a_0,\ldots,a_d \in K[z_l \st l \in J]$, then $a_0 = \ldots = a_d = 0$. If not, suppose (without loss of generality) that $a_0 \neq 0$. Then $v(a_0) = 0$, while $v(a_k (X - c_jY)^k) = v(a_k) + k = k > 0$ for each $k > 0$. Hence $\infty = v(0) = v(\sum_{n = 0}^d a_n (X-c_jY)^n) = 0$, a contradiction. \qed

In view of Lemma \cp1.4, for each $J \subseteq I, j \in J$, each element of $K[z_l \st l \in J][X - c_jY]$ has a unique presentation as a polynomial in $X - c_jY$. Thus we have a natural valuation on this ring, given by $v'(\sum_{n = 0}^d a_n (X - c_jY)^n) = \min(n \st a_n \neq 0)$, and we may form the completion $K[z_l \st l \in J][[X - c_jY]]$ of this ring with respect to $v'$.

\proclaim{Proposition \cp1.5} Let $J \subseteq I$ and $j \in J$. Then $D_J = K[z_l \st l \in J][[X - c_jY]]$, and $v$ is given on $D_J$ by $v(\sum_{n = 0}^\infty a_n (X-c_jY)^n) = \min(n \st a_n \neq 0)$. \endproclaim \demo{Proof} By Lemma \cp1.4, $v$ coincides with $v'$ (given in the paragraph preceding this proposition) on $K[z_l \st l \in J][X - c_jY]$, hence $K[z_l \st l \in J][[X - c_jY]]$ is the completion of $K[z_l \st l \in J][X - c_jY]$ with respect to $v$, and $v$ coincides with $v'$ on the completion. Note that $K[z_l \st l \in J][X - c_jY] = K[z_l \st l \in J][X,Y]$ (since $Y = z_j(X-c_jY)$, $X = (1 + c_jz_j)(X-c_jY)$), hence by the definition of $D_J$ we are done. \qed

\proclaim{Lemma \cp1.6} Let $J \subseteq I$. Then $K[z_j \st j \in J] = \sum_{j \in J} K[z_j]$. \endproclaim \demo{Proof} For each $i \neq j \in J$ we have $$z_i \cdot z_j = {Y^2 \over (X-c_iY)\cdot(X-c_jY)} = {1 \over c_i - c_j}\cdot z_i + {1 \over {c_j-c_i}} z_j.$$ The claim now follows by induction on $|I|$. \qed

\proclaim{Proposition \cp1.7} Let $J, J' \subseteq I$. Then for each $f \in D_{J \cup J'}$ there exist $f_1 \in D_J$ and $f_2 \in D_{J'}$ with $v(f_1),v(f_2) \geq v(f)$, such that $f = f_1 + f_2$. \endproclaim \demo{Proof} Replace $J$ with $J \smallsetminus (J \cap J')$ to assume that $J \cap J' = \emptyset$. Moreover, without loss of generality $J,J'$ are non-empty. Choose $j \in J,j' \in J'$ and denote $A_J = K[z_l \st l \in J], A_{J'} = K[z_l \st l \in J'], A = K[z_l \st l \in J \cup J']$. By Proposition \cp1.5, $D_J = A_J[[X - c_jY]]$, $D_{J'} = A_J[[X - c_{j'}Y]]$, $D_{J \cup J'} = A[[X - c_jY]]$. Let $f = \sum_{n = m}^\infty a_n (X - c_jY)^n \in D_{J \cup J'}$, with $a_m \neq 0$. Then $v(f) = m$ by Proposition \cp1.5. By Lemma \cp1.6, $A = A_J + A_{J'}$. For each $n \geq m$, let $b_n \in A_J, b_n' \in A_{J'}$ such that $a_n = b_n + b_n'$. Denote $f_1 = \sum_{n = m}^\infty b_n (X - c_jY)^n$, $f_2 = f - f_1 = \sum_{n = m}^\infty b_n' (X - c_jY)^n$. Then $f_1 \in D_J$ and $v(f_1) \geq m$. It remains to prove that $f_2 \in D_{J'}$ and that $v(f_2) \geq m$. This follows by the following equality: $$f_2 = \sum_{n = m}^\infty b_n' (X - c_jY)^n = \sum_{n = m}^\infty b_n' ((X - c_{j'}Y) + (c_{j'} - c_j)Y)^n$$ $$ = \sum_{n = m}^\infty (b_n' (1 + (c_{j'} - c_j)z_{j'})^n) (X - c_{j'}Y)^n.$$ \qed

The next lemma is a variant of [HaH10, Lemma 3.3], allowing non-principal ideals.

\proclaim{Lemma \cp1.8} Let $R \subseteq R_1,R_2 \subseteq R_0$ be Noetherian domains such that $R_0 = R_1 + R_2$. Let $w$ be a discrete valuation on $\Quot(R_0)$ such that $R$ is complete with respect to $w$ and $w(x) \geq 0$ for all $x \in R_0$. Let $\frp,\frp_1,\frp_2,\frp_0$ be the centers of $w$ in $R, R_1,R_2, R_0$, respectively. Suppose that $\frp R_0 = \frp_0$ and $R/\frp = R_1/\frp_1 \cap R_2 / \frp_2$ (inside $R_0/\frp_0$). Then $R_1 \cap R_2 = R$. \endproclaim \demo{Proof} First, note that $\frp_0 = \frp_1 + \frp_2$. Indeed, since $R$ is Noetherian, $\frp$ is generated by elements $x_1, \ldots, x_n \in R$, hence $\frp_0 = \frp R_0$ is also generated by $x_1,\ldots,x_n$. Suppose $\sum_{i=1}^n a_i x_i \in \frp_0$, $a_1,\ldots,a_n \in R_0$. For each $1 \leq i \leq n$, write $a_i = b_i + b'_i$ with $b_i \in R_1$ and $b'_i \in R_2$. Then $\sum a_i x_i = \sum b_i x_i + \sum b'_i x_i \in \frp_1 + \frp_2$, since $\frp \subseteq \frp_1,\frp_2$.

Denote $S = R_1 \cap R_2$ and let $\frq$ be the center of $w$ at $S$. Then the sequence $0 \to S \to R_1 \times R_2 \to R_0 \to 0$ is exact (where the second map is the diagonal map and the third map is substraction). This sequence induces an exact sequence $0 \to S/\frq \to (R_1/\frp_1) \times (R_2/\frp_2) \to R_0/\frp_0 \to 0$. Indeed, the only non-trivial part in showing this is to check that the kernel of the substraction map is contained in the image of the diagonal map. Suppose $(x_1 + \frp_1, x_2 + \frp_2) \in (R_1/\frp_1) \times (R_2/\frp_2)$ satisfies $x_1 - x_2 \in \frp_0$. Since $\frp_0 = \frp_1 + \frp_2$ we may choose $y_1 \in \frp_1, y_2 \in \frp_2$ such that $x_1 - y_1 = x_2 - y_2$. Then $(x_1 + \frp_1, x_2 + \frp_2) = (x_1 -y_1 + \frp_1, x_2 - y_2 + \frp_2)$ belongs to the image of the diagonal map. Thus the sequence is exact.

On the other hand, by our assumptions the sequence $0 \to R/\frp \to (R_1/\frp_1) \times (R_2/\frp_2) \to R_0/\frp_0 \to 0$ is also exact. Hence the natural map $R/\frp \to S/\frq$ is an isomorphism. In particular, $S = R + \frp S$. By induction we have $S = R + \frp^kS$ for each $k \in \bbN$. Thus $R$ is $w$-dense in $S$, hence the completion of $R$ with respect to $w$ contains $S$. But by our assumptions, $R$ is complete, hence $R = S$. \qed

\proclaim{Lemma \cp1.9} Suppose $a_0 + \sum_{i \in I}\sum_{n = 1}^{d_i} a_{i,n}z_i^n = 0$, where $d_i \in \bbN$ and $a_0, a_{i,n} \in K$ for each $i,n$. Then $a_0 = a_{i,n} = 0$ for each $i,n$. \endproclaim \demo{Proof} Suppose there exists $i \in I, n \in \bbN$ such that $a_{i,n} \neq 0$. Without loss of generality, $n = d_i$. Since $X - c_iY$ is a prime element of $K[X,Y]$, it defines a discrete valuation on $K(X,Y)$, which we denote by $w$. We have $w(Y) = w(Y - c_jX) = 0$ for each $j \neq i$ in $I$. Thus $w(a_0 + \sum_{j \neq i}\sum_{n = 1}^{d_j} a_{j,n}z_j^n) \geq 0$, while $w(\sum_{n = 1}^{d_i} a_{i,n}z_i^n) = -d_i$. Thus $w(0) = w(a_0 + \sum_{j \in I}\sum_{n = 1}^{d_j} a_{j,n}z_j^n) = -d_i$, a contradiction. \qed

\proclaim{Proposition \cp1.10} Suppose $J,J' \subseteq I$. Then $D_J \cap D_{J'} = D_{J \cap J'}$. \endproclaim \demo{Proof} Clearly, $D_{J \cap J'} \subseteq D_J \cap D_{J'}$. For the converse inclusion, we distinguish between two cases. First suppose that $J \cap J' \neq \emptyset$, and fix $j \in J \cap J'$. Then $D_J = K[z_k \st k \in J][[X - c_jY]]$ and $D_{J'} = K[z_k \st k \in J'][[X - c_jY]]$, hence $D_J \cap D_{J'} = (K[z_k \st k \in J] \cap K[z_k \st k \in J'])[[X - c_jY]]$. By Lemma \cp1.9 $K[z_k \st k \in J] \cap K[z_k \st k \in J'] = K[z_k \st k \in J \cap J']$, hence $y \in D_{J \cap J'}$.

Now suppose that $J \cap J' = \emptyset$ and denote $R = K[[X,Y]] = D_\emptyset, R_1 = D_J, R_2 = D_{J'}, R_0 = D_{J \cup {J'}}$. Since $v(f) \geq 0$ for each $f \in K[z_j \st j \in J \cup J'][X,Y]$, we also have $v(f) \geq 0$ for each $f$ in the completion $R_0$. The ring $R$ is complete with respect to $v$, and $R = R_1 + R_2$ by Proposition \cp1.7. Let $\frp,\frp_1,\frp_2,\frp_0$ be the centers of $v$ at $R, R_1,R_2, R_0$, respectively. Then $\frp$ is generated by $X,Y$ and $\frp_0$ is generated by $X - c_jY$ for any $j \in J$, by Proposition \cp1.5. It follows that $\frp R_0 = \frp_0$. In order to apply Lemma \cp1.8, it remains to check that $R_1/\frp_1 \cap R_2/\frp_2 = R/\frp$ in $R_0/\frp_0$. Indeed, we have $R_1/\frp_1 = K[z_j \st j \in J]$, $R_2/\frp_2 = K[z_j \st j \in J']$, $R_0/\frp_0 = K$. By Lemma \cp1.9, we are done. \qed

The rings $D_J$, $J \subseteq I$, that we have defined in an algebraic manner, have a natural rigid-geometric interpretation:

\rem{Remark \cp1.11} Let $J \subseteq I$, $j \in J$. Denote $t = X - c_jY$. By Proposition \cp1.5, $D_J = K[z_l \st l \in J][[t]]$ is the $t$-adic completion of $K[z_l \st l \in J][t]$, thus $D_J[t^{-1}]$ is the $t$-adic completion of $K[z_l \st l \in J][t,t^{-1}]$. We have $K[z_l \st l \in J][t,t^{-1}] \subseteq K((t))[z_l \st l \in J] \subseteq D_J[t^{-1}]$, hence $D_J[t^{-1}]$ is the $t$-adic completion of $A := K((t))[z_l \st l \in J]$. Denote $T = K[[t]], F = K((t))$, and let $v_t$ be the $t$-adic valuation on $F$. Then for an element $$f = f_0 + \sum_{k \in J} \sum_{n=1}^\infty f_{kn}z_k^n \leqno(\cf1)$$ in A, where $f_0, f_{kn} \in F$ are almost all zero, we have $v_t(f) = \min_{kn}\{v_t(f_0), v_t(f_{kn})\}$.

Denote $s = {X \over Y}$. Then $s$ is a free variable over $F$, and note that $z_j = {1 \over s - c_j}$ and $v_t(c_l - c_k) = 0$ for all distinct $l,k \in J$. By [HaJ98, Lemma 3.3] the completion $D_J[t^{-1}]$ of $A$ is the ring of holomorphic functions on the affinoid $U = \bbP \hefresh (\bigcup_{l \in J}B(c_l))$, where $\bbP$ is the projective $s$-line and $B(c_l)$ is a disc of radius $1$ with center $c_l$ for each $l \in J$ (cf. [FrP04, \S2.2]). Moreover, each element $f \in D_J[t^{-1}]$ can be uniquely presented as in (\cf1), where $f_0 \in F$ and $\{f_{ln}\}_{n=1}^\infty$ is a null sequence in $F$ (with respect to $v_t$) for each $l \in J$. Thus, $D_J$ is the ring of holomorphic functions on $U$ having no pole at $t$. Its elements are of the form (\cf1), where the coefficients are now in $T$ (and $\{f_{kn}\}_{n=1}^\infty$ is a null sequence for each $k \in J$). In particular, $T[z_k \st k \in J]$ is dense in $D_J$. \endrem

Note that the interpretation given in Remark \cp1.11 is dependent on the change of variable $t = X -c_jY$ for some $j \in J$. Thus, for disjoint sets $J, J' \subset I$ it is difficult to use this presentation to study the way $D_J$ and $D_{J'}$ fit together in $D_I$ (in order to do that, as in Proposition \cp1.10, we needed to view these rings as completions at the common valuation $v$). However, this interpretation is useful to study inner properties of the rings $D_J$, and consequently, of the rings $Q_i$. Most notably, these rings are fields.

\proclaim{Corollary \cp1.12} Let $J \subseteq I, j \in J$. \condition (a) For each $0 \neq g \in D_J$, $D_J[(X-c_jY)^{-1}] = K((X-c_jY))[z_k \st k \in J] + gD_J[(X-c_jY)^{-1}]$. \condition (b) For each $f \in D_J$ there exist $p \in K[[X-c_jY]][z_j], u \in D_J^\times$ such that $f = pu$. \condition (c) The ring $Q_j$ is a field. \endproclaim \demo{Proof} In the notation of Remark \cp1.11, each element $f \in D_J[t^{-1}]$ can be written in the form $u \cdot p$, with $u \in D_J[t^{-1}]^\times$ and $p \in F[z_j]$, by [HaJ98, Lemma 3.7]. If $f \in D_J$ then we can multiply $u$ and $p$ with a power of $t$ to assume that $u \in D_J^{\times}$ and $p \in K[[t]][z_j]$. This proves (b). Part (a) is given by [HaJ98, Corollary 3.8]. To prove (c), we first claim that $\Quot(D_J) = E D_J$. Indeed, take ${f \over g}\in \Quot(D_J)$ with $0 \neq f,g \in D_J$. By (b), we may multiply $f$ and $g$ with an element of $D_J^{\times}$ to assume that $f \in D_J$ and $g \in K[[t]][z_j]$. Since $K[[X,Y]]$ is complete with respect to $\langle X,Y \rangle$ and $t \in \langle X,Y \rangle$, $K[[X,Y]]$ is complete with respect to $t$, hence $T \subseteq K[[X,Y]]$ and $F \subseteq E$. Since $z_j \in E$, we get that $F[z_j] \subseteq E$. Thus ${f \over g} = f \cdot {1 \over g} \in D_J \cdot E$. Applying this argument for $I \hefresh \{j\}$ instead of $J$ we have $Q_j = \Quot(D_{I \hefresh \{j\}}) = E D_{I \hefresh \{j\}}$, hence  $\Quot(Q_j) = \Quot(E D_{I \smallsetminus \{j\}}) = E \Quot(D_{I\smallsetminus\{j\}}) = E D_{I\smallsetminus\{j\}} = Q_j$ is a field. \qed

\proclaim{Proposition \cp1.13} For each $i \in I$, $Q'_i \cap Q_i = E$. \endproclaim \demo{Proof} Since $Q_i' = \bigcap_{j \neq i} Q_j$, the assertion is that $\bigcap_{j \in I} Q_j = E$. Indeed, let $y \in \bigcap_{j \in I} Q_j$. For each $j \in J$ write $y = {f_j \over q_j}$ with $f_j \in D_{I \smallsetminus \{j\}}$, $q_j \in K[[X,Y]]$. Taking a common denominator we may assume that $q_j$ is independent of $j$, and denote $q = q_j$ (for all $j \in J$). It suffices to prove that $qy \in K[[X,Y]] \subseteq E$. But $qy = q_jy = f_j \in D_{I \smallsetminus \{j\}}$ for all $j \in I$, hence $qy \in \bigcap_{j \in I} D_{I \smallsetminus \{j\}} = D_\emptyset = K[[X,Y]]$, by Proposition \cp1.10. \qed

The proof of the following proposition is based on that of [HaJ98, Corollary 4.4]. (We cannot use [HaJ98, Corollary 4.4] as it is, since condition (e') of that claim does not hold for $D_I$ itself.)

\proclaim{Proposition \cp1.14} Let $i \in I$, $n \in \bbN$ and let $b \in \GL_n(Q)$. There exist $b_1 \in \GL_n(Q_i)$ and $b_2 \in \GL_n(Q'_i)$ such that $b = b_1 \cdot b_2$. \endproclaim \demo{Proof} Denote by $|\cdot|$ the absolute value on $Q$ that corresponds to $v$. The rings $A = D_I$, $A_1 = D_{I \smallsetminus \{i\}}, A_2 = D_{\{i\}}$ are complete with respect to $|\cdot|$ and Proposition \cp1.7 asserts that condition (d') of [HaJ98, Example 4.3] holds for these rings. We extend $|\cdot|$ to the maximum norm $||\cdot||$ on $M_n(Q)$, as in [HaJ98, Example 4.3]. Then $M_n(A), M_n(A_1), M_n(A_2)$ are complete with respect to $||\cdot||$ and condition (d) of [HaJ98, \S4] holds. By Cartan's Lemma [HaJ98, Lemma 4.2], for each $a \in \GL_n(A)$ with $||a - 1|| < 1$ there exist $a_1 \in \GL_1(A_1)$, $a_2 \in \GL_1(A_2)$ such that $a = a_1 \cdot a_2$.

Denote $E_1 = \Quot(A_1) = Q_i, E_2 = \Quot(A_2) = Q_i'$. In order to factor $b$ (which need not be in $\GL_n(A)$), let $t = X - c_iY$, $T = k[[t]]$. By Remark \cp1.11 (for $J = I)$ $A_0 = T[z_k \st k\in I]$ is a dense subring of $A$, and by Corollary \cp1.12(b) there exists $h \in A_0$ such that $hb \in \M_n(A)$. If $hb = b_1b_2'$ with $b_1 \in \GL_n(E_1)$ and $b'_2 \in \GL_n(E_2)$, then $b = b_1b_2$ with $b_2 = {1 \over h} b_2' \in \GL_n(E_2)$. So we may assume that $b \in \M_n(A)$.

Let $0 \neq d = \det(b) \in A$. By Corollary \cp1.12(b) there are $0 \neq g \in A_0$ and $u \in A^\times$ such that $d = gu$. Let $b'' \in \M_n(A)$ be the adjoint matrix of $b$, so that $bb'' = d1$. Let $b' = u^{-1}b''$. Then $b' \in \M_n(A)$ and $bb' = g1$. Put $$V = \{a' \in \M_n(A[t^{-1}]) \st ba' \in gM_n(A[t^{-1}])\}, V_0 = V \cap M_n(A_0[t^{-1}]).$$ Then $V$ is an additive subgroup of $M_n(A[t^{-1}])$ and $gM_n(A[t^{-1}]) \leq V$. By Corollary \cp1.12(a) $M_n(A[t^{-1}]) = M_n(A_0[t^{-1}]) + gM_n(A[t^{-1}])$, hence $V = V_0 + gM_n(A[t^{-1}])$. Since $A_0$ is dense in $A$, $gM_n(A_0[t^{-1}])$ is dense in $gM_n(A[t^{-1}])$. It follows that $V_0 = V_0 + gM_n(A_0[t^{-1}])$ is dense in $V = V_0 + gM_n(A[t^{-1}])$. As $b' \in V$, there exists $a_0 \in V_0$ such that $||b' - a_0|| < {|g| \over ||b||}$.
Put $a = {1 \over g} a_0 \in M_n(Q)$. Then $ba \in M_n(A[t^{-1}])$ and $||1 - ba|| = ||{1 \over g}b(b'-a_0)|| \leq {1 \over |g|}||b||\cdot||b'-a_0|| < 1$. Hence $||ba|| = 1$, so each entry in $ba$ has a non-negative value at $v$. By Remark \cp1.11, $v$ coincides with the $t$-adic valuation on $A$, hence all the entries of $ba$ belong to $A$. Thus $ba \in M_n(A)$, and since $||1 - ba|| < 1$ and $M_n(A)$ is complete, $ba \in \GL_n(A)$. In particular, $\det(a) \neq 0$ and hence $a \in \GL_n(\Quot(A_0)) \subseteq \GL_n(E_2)$. By the first paragraph, there exist $b_1 \in \GL_n(A_1) \subseteq \GL_n(E_1), b'_2 \in \GL_n(A_2)$ such that $ba = b_1 b_2'$. Then $b_2 = b_2'a^{-1} \in \GL_n(E_2)$ satisfies $b = b_1b_2$. \qed

\proclaim{Lemma \cp1.15} Let $J \subseteq I, j \in J$, $t = X - c_jY$. \condition (a) Suppose $p = \sum_{l=0}^d p_l z_j^l \in K[[t]][z_j]$ is a polynomial such that $v(p_1) =  0$ and $v(p_l) > 0$ for each $l > 1$. Then $p$ is prime in $D_J[t^{-1}]$. \condition (b) The ring $D_J[t^{-1}]$ is a unique factorization domain. \condition (c) For each $a,b,c \in K^\times$ with $a \neq -b$ and $2 \leq m \in \bbN$, the elements $1 + a z_j + t^{m-1}z_j^m$, $1 + b z_j - t^{m-1}z_j^m$, $1 + c z_j$ are non-associate primes of $D_J[t^{-1}]$. \endproclaim \demo{Proof} Denote $F = K((t))$. Then $D_J[t^{-1}] = F\{z_k \st k \in J\}$ (Remark \cp1.11). Viewing $p$ as an element of $F\{z_j\}$, it is regular of pseudo degree $1$ [HaV96, Definition 1.4], hence by [HaV96, Corollary 1.7] we have $p = u \cdot q$, where $u \in F\{z_j\}^\times \subseteq D_J[t^{-1}]^\times$ and $q = q_0 + z_j\in F[z_j]$ is a linear polynomial with $v(q_0) \geq 0$. Thus to prove (a), it suffices to show that $q$ is prime in $D_J[t^{-1}]$. Without loss of generality $q_0 \neq 0$, and we denote $c = c_j - {1 \over q_0}$. Then $q = z_j - {1 \over c - c_j}$, hence by [Pa08, Lemma 6.4(a)] (with $D$, $r$, $1$ there replaced by $F$, $1$, $j$ here) $q$ generates the kernel of an epimorphism from $D_J[t^{-1}]$ onto a domain (acutally a field here), hence $q$ is prime. This proves (a).

Since $D_J[t^{-1}]$ is a principal ideal domain [HaJ98, Proposition 3.9], part (b) follows.

By part (a), $r = 1 + a z_j + t^{m-1}z_j^m$, $r' = 1 + b z_j - t^{m-1}z_j^m$, $s = 1 + c z_j$ are primes of $D_J[t^{-1}]$. If $s | r$, then $-{1 \over c}$ is a root of $r$, a contradiction. Thus $r$, $s$ (and similarly, $r'$, $s$) are non-associates.

If $r | r'$ then $r | r + r'$. By the argument of the preceding paragraph, $r + r' = 2 + (a+b)z_j$ is a prime, non associate to $r$, a contradiction. This proves (c). \qed


\proclaim{Lemma \cp1.16} Let $K$ be a field that contains a primitive $q$-th root of unity, for some $q \in \bbN$. Let $v$ be a discrete valuation on $K$ which is trivial on the prime field of $K$, and let $a \in K$ with $v(a) = 0$. Suppose $L = K(a^{1 \over q})$ is a Kummer extension of $K$, and that $L/K$ is unramified at $v$. Then $v(x^\sig) = v(x)$ for each $x \in L$ and $\sig \in \Gal(L/K)$. \endproclaim \demo{Proof} Extend $v$ arbitrarily to $L$, let $O$ be the valuation ring of $v$ in $K$, and $O'$ the valuation ring of $v$ in $L$. Since $K$ contains a primitive $q$-th root of unity, $q$ is not divisible by $p = \chr(K)$. Thus $d = \disc(T^q - a,K) = ka^{q-1}$, where $k \in \bbZ$ is not divisible by $p$. Hence $v(d) = 0$, and by [FrJ05, Lemma 6.1.2] we have $O' = O[a^{1 \over q}]$. Put $\alp = a^{1 \over q}$ and let $x = \sum_{i=0}^{q-1}b_i \alp^n \in K$, with $b_0,\ldots,b_{q-1} \in K$. We claim that $v(x) = \min_i(v(b_i))$. Indeed, since $L/K$ is unramified at $v$, we may multiply $x$ by a power of a uniformizer of $v$ in $K$, to assume that $v(x) = 0$. Since $O' = O[\alp]$, $v(b_i) \geq 0$ for each $0 \leq i \leq {q-1}$. On the other hand $v(x) \geq \min_i (v(b_i \alp^i)) = \min_i(v(b_i))$, since $v(\alp) = {1 \over n} v(a) = 0$. Thus $v(b_i) = 0$ for some $0 \leq i \leq {q-1}$, hence $v(x) = \min_i(v(b_i))$.

Now, let $\sig \in \Gal(L/K)$ and let $x = \sum_{i=0}^{q-1}b_i \alp^n \in K$, with $b_0,\ldots,b_{q-1} \in K$, be an arbitrary element. We have $\alp^\sig = \zeta \alp$, where $\zeta$ is some $q$-th root of unity. Then $v(x^\sig) = v(\sum_{i=0}^{q-1} b_i \zeta^i \alp^i) = \min_i(v(b_i\zeta^i)) = \min_i(v(b_i)) = v(x)$. \qed

Recall that given a field $K$, any $K$-central simple algebra $A$ is of the form $\Mat_n(D)$ for some $K$-division algebra $D$. The index of $A$ is defined to be $\ind(A)=\sqrt{\dim_KD}$. So, $A$ is a division algebra if and only if $\ind(A)=\sqrt{\dim_K A}$. Let us denote Brauer equivalence by $\sim$ and the exponent of an algebra (its order in the Brauer group) by $\exp$. A subfield $F$ of $D$ is a maximal subfield of $D$ if and only if $\ind(D)=[F:K]$. Furthermore, a field $F$ is a maximal subfield of $D$ if and only if $\ind(D)=[F:K]$ and $F$ splits $D$.

The proof of the next proposition is partially based on that of [HHK10, Proposition 4.3].

\proclaim{Proposition \cp1.17} Fix $i \in I$, and let $H$ be an abelian $p$-group of rank at most $2$, where $p \neq \chr(K)$. Suppose $K$ contains an $|H|$-th primitive root of unity. Let $E'$ be a finite extension of $E$. Then there exists an $H$-Galois extension $F_i/E$ such that:
\condition (a) $F_i \subseteq Q_i'$.
\condition (b) $F_i$ is contained as a maximal subfield in an $E$-division algebra $D'_i$, and $D'_i \tensor_E E'Q_i$ remains a division algebra (where $E'Q_i$ is the compositum of $E'$ and $Q_i$ in the algebraic closure of $Q$).
\endproclaim \demo{Proof} Write $H = C_q \times C_q'$, where $q,q'$ are powers of $p$. For each $k \in \bbN$, the elements $X - c_iY + Y^k$, $X + c_iY - Y^k$ are irreducible and hence prime in the unique factorization domain $K[[X,Y]]$. Only finitely many primes of $K[[X,Y]]$ are ramified at $E'/E$, hence for a sufficiently large $2 \leq k \in \bbN$, $f = X - c_iY + Y^k$ and $g = X + c_iY - Y^k$ are unramified at $E'/E$. That is, the corresponding valuations $v_f,v_g$ are unramified. Let $a = {f \over X - c_iY}$, $b = {g \over X - c_iY}$. Clearly $v_f(X - c_iY) = v_f(g) = 0$, hence $v_f(a) = 1, v_f(b) = 0$. Similarly, $v_g(a) = 0, v_g(b) = 1$. Consider the polynomial $h(U) = U^q - a$ over $D_{\{i\}} = K[z_i][[X-c_iY]]$. Note that $a = 1 + z_i^k(X - c_iY)^{k-1}$, hence $h(1) \in \ (X-c_iY)D_{\{i\}}$ and $h'(1) = q \in K^\times \subseteq D_{\{i\}}^\times$. By the ring version of Hensel's Lemma (for the ideal $(X-c_iY)D_{\{i\}}$) $h(U)$ has a root $s \in D_{\{i\}}$. Note that $v_f(s) = {1 \over q} \notin \bbZ$, hence $s \notin E$. Since $K$ contains a primitive $|H|$-th root of unity, it contains a primitive $q$-th root of unity. By Kummer theory $E(s)/E$ is a Galois extension with group $C_q$. Similarly, there exists $s' \in D_{\{i\}}$ satisfying $(s')^{q'} = b$, and $E(s')/E$ is Galois with group $C_{q'}$. Let $F_i = E(s,s') \subseteq Q_i'$.

Since $v_f(a) = 1$, $h(U)$ is irreducible over $E$, by Eisenstein's criterion. Denoting the reduction modulo $g$ by $\bar{\cdot}$, $\bar{h}(U) = U^q - \bar{a}$ is separable, since $\bar{a} \neq 0$. Thus by [FrJ05, Lemma 2.3.4], $E(s)/E$ is unramified at $v_g$. Clearly, $E(s')/E$ is totally ramified at $v_g$. Thus $E(s),E(s')$ are linearly disjoint over $E$, hence $\Gal(F_i/E) = H$.

Let $D_i'$ be the quaternion algebra $(a,b)_{qq'}$. Note that $D_i'$ can be also viewed as the cyclic algebra $(E({a}^{1\over qq'})/E,\sigma,b)$, for some generator $\sigma$ of $\Gal(E({a}^{1\over qq'})/E)$. We claim that $F_i$ is a maximal subfield of $D_i'$. As $\dim_E(D_i')=[F_i:E]^2$ we only need to show that $F_i$ splits $D_i'$. By [Rei75, Theorem 30.8],  $D_i'\otimes_E E(s)\sim (E({s}^{1\over q'})/E(s),\sigma^q,b)$ and thus $D_i'\otimes_E F_i\sim  (F_i({s}^{1\over q'})/F_i,\sigma^q,b)$. By [Rei75, Theorem 30.4], the latter algebra splits if and only if $b\in N_{F_i({s}^{1\over q'})/F_i}(F_i({s}^{1\over q'}))$. This holds since $b=N_{F_i({s}^{1\over q'})/F_i}(s')$ which shows that $F_i$ is indeed a maximal subfield of $D_i'$.

Choose $j \in I \hefresh \{i\}$, denote $t = X - c_jY$ and $r = 1 + (c_j + c_i)z_j -t^{k-1}z_j^k, r' = 1 + (c_j - c_i)z_j + t^{k-1}z_j^k$. By Lemma \cp1.15(c) $r,r'$ are non-associate prime elements in $D_{I\hefresh\{i\}}[t^{-1}]$, so they define discrete valuations $v_r, v_{r'}$ on $Q_i = \Quot(D_{I\hefresh\{i\}}) = \Quot(D_{I\hefresh\{i\}}[t^{-1}])$ such that $v_r(r') = v_{r'}(r) = 0$. By Lemma \cp1.15(c) we also have $v_{r'}(1 + (c_j - c_i)z_j) = v_r(1 + (c_j - c_i)z_j) = 0$.

Note that $$b = {X  - c_jY + (c_j+c_i)Y -Y^k \over X - c_jY + (c_j-c_i)Y} = {t + (c_j + c_i)tz_j - t^jz_j^k \over t + (c_j - c_i)tz_j} = {r \over 1 + (c_j - c_i)z_j}.$$ Similarly, $a = {r' \over 1 + (c_j - c_i)z_j}$. Thus $v_r(b) = 1, v_{r'}(b) = 0, v_{r}(a) = 0, v_{r'}(a) = 1$. Then the polynomial $U^{qq'} - a$ is irreducible over $D_{I\hefresh\{i\}}$, by Eisenstein's Criterion (using $v_{r'}$). Thus $Q_i({a}^{1\over qq'})/Q_i$ is unramified at $v_r$ (again by [FrJ05, Lemma 2.3.4]), hence so is $E'Q_i({a}^{1\over qq'})/E'Q_i$.

Only finitely many primes of the unique factorization domain $D_{I\hefresh\{i\}}[t^{-1}]$ (Lemma \cp1.15(b)) are ramified at the finite extension $E'Q_i/Q_i$, hence without loss of generality, we may assume that $E'Q_i/Q_i$ is unramified at $v_{r'}$ (by possibly choosing an even larger $k$ before hand). On the other hand, $Q_i(a^{1\over qq'})/Q_i$ is totally ramified at $v_{r'}$, hence $[E'Q_i(a^{1\over qq'}):E'Q_i] = [Q_i(a^{1\over qq'}):Q_i] = qq'$.

We can now show that $D_i'\otimes_E E'Q_i$ is a division algebra. A sufficient condition for this to hold is that $\exp(D_i'\otimes_E E'Q_i)=qq'$. This happens if and only if for every $1\leq m\leq qq'-1$ the algebra $(E'Q_i(a^{1\over qq'})/E'Q_i,\sigma,b^m)\sim (D_i'\otimes_E E'Q_i)^m$ does not split. Let $N$ denote the norm $N_{E'Q_i(a^{1\over qq'})/E'Q_i}$. The splitting of the latter cyclic algebra is equivalent to having $b^m\not\in N(E'Q_i(a^{1\over qq'})^\times)$ for all $1\leq m\leq qq'-1$.

As $E'Q_i({a}^{1\over qq'})/E'Q_i$ is unramified at $v_r$, we have $v_r(x) = v_r(x^\sig)$ for each $x \in E'Q_i({a}^{1\over qq'})$, by Lemma \cp1.16. Hence $v_r(N(x))=\sum_{l=0}^{qq'-1}v_r(x^{\sigma^l}) = qq'v_r(x)$ for all $x \in E'Q_i({a}^{1\over qq'})$. Since $v_r(b)=1$, $b^m\not\in N(E'Q_i({a}^{1\over qq'})^\times)$ for all $1\leq m\leq qq'-1$. \qed

%
%
%
%

\head \cs2. Patching and admissibility \endhead

We have established the patching machinery needed to prove our main theorem (Theorem \cp2.8). We first recall some general properties of induced algebras and Frobenius algebras.

\reminfo{Remark \cp2.1}{Induced Algebras} Let $G$ be a finite group and $H\leq G$. Let $P/Q$ be a finite field extension with $H=\Gal(P/Q)$. 
Let $N=\Ind^G_{H} P = \Big\{ \sum_{\tet\in G} a_\tet \tet \st a_\tet \in P, \ a_\tet^\tau=a_{\tet\tau} \hbox{ for all } \tet\in G, \ \tau \in H \Big\}$ be a ring with respect to point-wise addition and multiplication. 
Then $P$ can be embedded as a subring of $N$ by choosing representatives
$\tet_1=1,...,\tet_k$ of $H\setminus G$ and sending an element $x\in P$ to  $\sum_{i=1,\tau\in H}^k x^\tau
{\tet_i}\tau$. Furthermore, by choosing different representatives $N$ can be presented as a direct some of copies of $P$. If $P$ splits a central simple $Q$-algebra $A$ then $\Ind_{H}^GP\otimes_Q A\cong \Ind_{H}^GP\otimes_P (P\otimes_Q A)\sim \Ind_{H}^GP\otimes_P P$ which is a split separable (Azumaya) algebra over $\Ind_H^GP$ (for a definition of a split separable algebra over a ring see [DI71, \S5]).  \endrem

The following definition, remark and proposition will be useful in the sequel and all appear in [Jac96, Section 2.1].

\reminfo{Definition \cp2.2}{Frobenius algebras} Let $F$ be a field. An $F$-algebra $A$ is a Frobenius algebra if $A$ contains a hyperplane $H$ that contains no non-zero one sided ideals of $A$.

\rem{Remark \cp2.3} An algebra $A_1\oplus \ldots \oplus A_s$ is Frobenius if and only if $A_i$ is Frobenius for each $1 \leq i \leq s$. Any algebra $F[a]$ (with a single generator) is Frobenius. Let $L/K$ be an $H$-extension. By Remark \cp2.1, $\Ind_H^GL$ can be decomposed into a sum of copies of $L$ and it follows that $\Ind_H^GL$ is a Frobenius algebra. \endrem


\proclaim{Proposition \cp2.4 [Jac96, Theorem 2.2.3]} Let $A$ be an $F$-central simple algebra and $K$ a commutative Frobenius subalgebra of $A$ such that $\dim_F A=[K:F]^2$. Then any embedding of $K$ into $A$ can be extended to an inner automorphism of $A$. \endproclaim

\proclaim{Lemma \cp2.5} Let $R$ be a complete local domain of dimension $r$, with a separably closed residue field $C$. Then $R$ is a finite module over a subring of the form $C[[X_1,...,X_r]]$. \endproclaim \demo{Proof} By Cohen's structure theorem [Mat86, \S29], $R$ is finitely generated over a subring of the form $B = K[[X_1,\ldots,X_n]]$ or of the form $B = \bbZ_p[[X_1,\ldots,X_n]]$. Then $C$ is a finite extension of the residue field of $B$. In particular, the residue field of $B$ is infinite, disqualifying the option $B = \bbZ_p[[X_1,\ldots,X_n]]$. Since $\dim B = \dim R = r$, we have $B = K[[X_1,...,X_r]]$, and $C$ is the separable closure of $K$.

We claim that $R$ contains $C$. To prove this, it suffices to show that every separable polynomial $q(Z) \in K[Z]$ has a root in $R$. Indeed, since $C$ is separably closed, $q(Z)$ has a root $\zeta \in C$, and $q'(\zeta) \neq 0$. Let $\fr{m}$ be the maximal ideal of $R$, and denote the reduction modulo $\fr{m}$ by $\bar{\cdot}$. Take $a \in R$ such that $\bar{a} = \zeta$. Then $a \notin \fr{m}$, hence $a \in R^\times$. Thus $q(a) \in \fr{m}$ and $q'(\bar{a}) = q'(\zeta) \neq 0$, hence $q'(a) \in R \hefresh \fr{m} = R^\times$. By Hensel's Lemma, $q(Z)$ has a root in $R$.

Since $[C:K]$ is finite, $C(K[[X_1,...,X_r]]) = C[[X_1,...,X_r]]$ is contained in $R$, and $R$ is finite over $C[[X_1,...,X_r]]$ (since it is finite over $K[[X_1,...,X_r]]$). \qed

The final ingredient we need in order to prove our main theorem is patching of central simple algebras. The content of the following proposition is essentially given in [HHHK10], but for specific fields $Q_i$, while here we present it for general fields satisfying a matrix factorization property. We note that [HHK10, Theorem 4.1] uses the terminology of categories and equivalence of categories. Here we prefer a more explicit presentation, working with vector spaces and bases, as in [HaJ98]. The proof of the proposition combines the proof of [HaJ98, Lemma 1.2] (where a more restricted assertion is made for specific types of algebras), and the proof of [HaH09, Theorem 7.1(vi)] (where the assertion is made for specific types of fields).

\proclaim{Proposition \cp2.6} Let $I$ be a finite set. For each $i \in I$ let $Q_i$ be a field contained in a $Q_i$-algebra $A_i$. Let $Q$ be a field containing $Q_i$ for each $i \in I$, and contained in a $Q$-algebra $A_Q$ which contains $A_i$ for each $i \in I$. Moreover, suppose that $A_iQ = A_Q$ and $\dim_{Q_i}A_i = \dim_{Q} A_Q$ for each $i \in I$. Finally, suppose that:

\condition (\cf2) For each $B \in \GL_n(Q)$ there exist $B_i \in \GL_n(Q_i)$ and $B_i' \in \GL_n(\bigcap_{j \neq i}Q_j)$ such that $B = B_i B_i'$.

Then, denoting $E = \bigcap_{i \in I} Q_i$, $A = \bigcap_{i \in I}A_i$ is an $E$-algebra satisfying $AQ_i = A_i$ for each $i \in I$. Moreover, if each $A_i$ is central simple, then so is $A$. \endproclaim \demo{Proof} For each $i \in I$, let $\calC_i$ be a basis for $A_i$ over $Q_i$. Since $A_iQ = A_Q$, $\Span_Q(\calC_i) = A_Q$, and since $\dim_{Q_i} A_i = \dim_Q A_Q$, $\calC_i$ is a basis for $A_Q$ over $Q$, for each $i \in I$. We now construct a basis $\calC$ for $A_Q$ over $Q$, which is also a basis for $A_i$ over $Q_i$, for {\bf all} $i \in I$.

For each subset $J$ of $I$, we find by induction on $|J|$, a basis $\calV_J$ for $A_Q$ over $Q$ which is also a basis for $A_j$ over $Q_j$, for each $j \in J$. Then for $I = J$ we will get the basis $\calC$.

If $J = \emptyset$ there is nothing to prove. Suppose that $|J| \geq 1$, choose $k \in J$ and let $J' = J \hefresh \{k\}$. By assumption there is a basis $\calV_J$ for $A_Q$ over $Q$ which is a basis for $A_i$ over $Q_i$ for each $i \in J'$. Since $\calC_i$ is a common basis for $A_Q$ and $A_i$, there is a matrix $B \in \GL_n(Q)$ such that $\calC_k B = \calV_{J'}$. By (\cf2), there exist $B_k \in \GL_n(Q_k)$ and $M \in \GL_n(\bigcap_{k \neq j \in I} Q_j) \subseteq \bigcap_{j \in J'} \GL_n(Q_j)$ such that $B = B_kM$. Put $\calV_J = \calV_{J'}M^{-1}$. Then $\calV_J$ is a basis for $A_Q$ over $Q$ which is also a basis for $A_j$ over $Q_j$ for each $j \in J'$. Moreover, $\calV_J$ is also a basis for $A_k$ over $Q_k$, since $\calV_J = \calV_kBM^{-1} = \calV_kB_k$. This completes the induction.

The existence of the common basis $C$ implies that $AQ_i = A_i$ for each $i \in I$. Now suppose that $A_i$ is central simple for each $i \in I$. Note that $E$ is contained in the center of $A$, $Z(A)Q_i$ is contained in the center of $A_i$ and hence $\dim_EZ(A) \leq \dim_{Q_i}(Z(A)Q_i)=1$ which shows that $Z(A)=E$. Thus $A$ is central over $E$. Similar dimension reasoning shows that a non-trivial ideal $I$ of $A$ extends to a non-trivial ideal of $A Q_i = A_i$. Thus $A$ is also simple. \qed

\proclaim{Proposition \cp2.7} Let $R$ be a complete local domain of dimension $2$, with a separably closed residue field. Let $G$ be a finite group of order not divisible by $\chr(K)$, whose Sylow subgroups are abelian of rank at most 2. Then $G$ is admissible over $\Quot(R)$. \endproclaim \demo{Proof} By Lemma \cp2.5, $R$ is a finite module over a subring of the form $B = K[[X,Y]]$, where $K$ is separably closed. Denote $E = \Quot(B) = K((X,Y))$, $E' = \Quot(R)$. Then $E'$ is a finite extension of $E$.

\subdemoinfo{Construction A}{A patching data} Let ${(p_i)}_{i \in I}$ be the prime factors of $n = |G|$, for some index set $I$. For each $i \in I$, let $G_i$ be a $p_i$-Sylow subgroup of $G$. Apply Construction \cp1.1 to obtain rings $Q_i$, $i \in I$, contained in the common field $Q$. Since $K$ is separably closed, it contains all roots of unity of order not divisible by $p = \chr(K)$. Thus for each $i \in I$, we may apply Proposition \cp1.17 to obtain a Galois extension $F_i/E$ with group $G_i$, such that $F_i \subseteq Q_i'$ and $F_i$ is contained as a maximal subfield in a division $E$-algebra $D_i'$. Moreover, $D'_i \tensor_E E'Q_i$ remains a division algebra. Thus $D_i := D'_i \tensor_E Q_i$ is also a division algebra. By Corollary \cp1.12(c), $Q_i$ is a field for each $i \in I$. Put $P_i = F_iQ_i$. Since $F_i$ splits $D'_i$, $P_i$ splits $D_i$, and since $[P_i:Q_i]=[F_i:E]=\ind(D'_i)=\ind(D_i)$, $P_i$ is a maximal subfield of $D_i$.

Let $\calE = (E, F_i, Q_i, Q;G_i, G)_{i \in I}$. Then $\calE$ is a patching data [HaJ98, Definition 1.1]. Indeed, conditions (2a), (2b) and (2d) of [HaJ98, Definition 1.1] have been established. Conditions (2c) and (2e) are given by Proposition \cp1.13 and Proposition \cp1.14, respectively.

\subdemoinfo{Part B}{Induced Algebras [HaJ98, \S1]} Consider the induced algebra $N = \Ind^G_1 Q$ of dimension $n$ over $Q$, and the $Q_i$-subalgebra $N_i = \Ind ^G_{G_i} P_i$ for each $i \in I$ (Remark \cp2.1). Then $G$ acts on $N$
by $\Big(\sum_{\tet\in G}a_\tet\tet\Big)^\sig = \sum_{\tet\in G}a_\tet\sig^{-1}\tet=\sum_{\tet\in G}a_{\sig\tet}\tet$ for each $\sig\in G$. The field $Q$ is embedded diagonally in $N$, which induces an embedding of $Q_i$ in $N_i$, for each $i \in I$. We view these embeddings as containments. By [HaJ98, Lemma 1.2] there is a basis for $N$ over $Q$, which is also a basis for $N_i$ over $Q_i$, for each $i \in I$. In particular, we have $N_iQ = N$ for each $i \in I$. By [Ha98, Lemma 1.3] $F=\inter_{i\in I} N_i$ is a Galois field extension of $E$ with group $G$, and there exists an $E$-embedding of $F$ into $Q$. Denote the image of $F$ under this embedding by $F'$.

%
%
%
%
%

\subdemoinfo{Part C}{Division algebras} It remains to prove that the extension $F'/E$ is adequate. Let $A_Q = \Mat_n(Q)$, and for each $i \in I$ let $n_i = [G:G_i]$. As $A_Q$ is split of dimension $n^2$ and $N$ is of dimension $n$ over $Q$, we also have an embedding of $N$ into $A_Q$. We view $N$ as a subalgebra of $A_Q$ via this embedding.

Fix $i \in I$. Since $P_i = F_iQ_i$ splits $D_i$, it follows by Remark \cp2.1 that $N_i$ also splits $D_i$. Moreover, by [DI71, Theorem 5.5] there is a central simple $Q_i$-algebra $A_i$ which is Brauer equivalent to $D_i$, in which $N_i$ embeds as a maximal commutative separable $Q_i$-subalgebra so that $\dim_{Q_i}(A_i)=\dim_{Q_i}(N_i)^2=n^2$. We view $N_i$ as contained in $A_i$ via this embedding.

Since $P_i$ splits $D_i$, we have $$D_i \otimes_{Q_i}Q \cong (D_i \otimes_{Q_i} P_i) \otimes_{P_i} Q \cong \Mat_{n \over n_i}(P_i) \otimes_{P_i} Q \cong \Mat_{n \over n_i}(Q).$$ Since $\ind(D_i)= {n \over n_i}$ and $\dim_Q A_Q = n^2$ we get that $A_i \cong M_{n_i}(D_i)$. Thus we have $A_i \otimes_{Q_i}Q \cong M_{n_i}(D_i)\otimes_{Q_i}Q \cong M_n(Q) = A_Q$, and we denote the induced $Q$-algebras isomorphism $A_i \otimes_{Q_i}Q \to A_Q$ by $\psi_i$. We cannot identify these two algebras via this isomorphism, since it might not be compatible with the containments $N_i \subseteq A_i$ and $N \subseteq A_Q$. This compatibility problem can be settled similarly to [HHK10, Lemma 4.2]:

By Part B we have $N = N_iQ$ and $\dim_{Q_i}N_i = \dim_QN = n$. Thus we have an isomorphism $\del_i \colon N =N_iQ \to N_i \otimes_{Q_i}Q$ for which the following diagram commutes:

$$ \xymatrix@=15pt{
A_i \ar@{-}[dd]\ar@{-}[rrrr] &&&& A_i \otimes_{Q_i}Q \\
\\
N_i \ar@{-}[rr] && N \ar[rr]^{\del_i} && N_i \otimes_{Q_i}Q \ar[uu]^{\id \otimes_{Q_i}Q}}
\leqno(\cf3) $$

By Remark \cp2.3, $N = \Ind_1^G Q$ is a Frobenius (commutative) subalgebra of $A_Q$. By Proposition \cp2.4, the embedding $\psi_i(\id \otimes_{Q_i}Q)\del_i \colon N \to A_Q$ extends to an inner automorphism $\alp_i$ of $A_Q$. Denote $\psi_i' = \alp_i^{-1}\psi_i$. Then $\alp_i^{-1}\psi_i(\id \otimes_{Q_i}Q)\del_i$ is the identity map on $N = N_iQ$, so we have the following commutative diagram:

$$ \xymatrix@=15pt{
A_i \otimes_{Q_i}Q \ar[rr]^{\psi_i'} && A_Q \\
\\
N_iQ \ar@{=}[rr] \ar[uu]^{(\id \otimes_{Q_i}Q)\del_i} && N \ar[uu]}
\leqno(\cf4) $$

Combining diagrams (\cf3) and (\cf4), we get the following commutative diagram

$$ \xymatrix@=15pt{
A_i \ar@{-}[rr] \ar@{-}[dd] && A_i \otimes_{Q_i}Q \ar[rr]^{\psi_i'} && A_Q \\
\\
N_i \ar@{-}[rr] && N_iQ \ar@{=}[rr] \ar[uu]^{(\id \otimes_{Q_i}Q)\del_i} && N \ar[uu]}
$$

This diagram gives an embedding $A_i \to A_Q$ which is compatible with the containments $N_i \subseteq A_i$ and $N \subseteq A_Q$, so we may now identify $A_i$ as a subring of $A_Q$, via this embedding. Moreover, since $\psi_i'$ is an isomorphism, we have $A_i \otimes_{Q_i} Q = A_iQ = A_Q$ by this identification. The following diagram explains the containment relations:

%
%

$$ \xymatrix@=15pt{
&& A_i \ar@{-}[dll] \ar@{-}[rr] \ar@{-}[dd] && A_Q \ar@{-}[dl]\\
N_i \ar@{-}[dd] \ar@{-}[rrr] &&& N \ar@{-}[dd]\\
&& D_i \ar@{-}[dl] \ar@{-}[dd] \\
Q_i \ar@{-}[dd] \ar@{-}[r] & P_i \ar@{-}[dd] \ar@{-}[rr] && Q \ar@{-}[dd] \\
&& D_i' \ar@{-}[dl] \\
E \ar@{-}[r] & F_i \ar@{-}[rr] && Q'_i }
$$

Let $A = \inter_{i \in I} A_i$. By Proposition \cp2.6, $A$ is a central simple $E$-algebra for which $AQ_i = A_i$ for each $i \in I$. In particular, $A =M_k(D)$ for some division algebra $D$ of index ${n\over k}$. 

Now, $D\otimes_E{Q_i}$ is Brauer equivalent to $A\otimes_E Q_i\cong A_i$ which is Brauer equivalent to $D_i$. Thus, ${n\over n_i}=\ind(D_i)|\ind(D)$ for each $i \in I$ and $n=\lcm_{i}({n\over n_i})|\ind(D)$. It follows that $k=1$ and $A$ is a division algebra. Naturally, $F$ is a subfield of $A$ and $\ind(A)=[F:E]$. It follows that $F$ is a maximal subfield of the division algebra $A$.

By choosing a basis for $A/F$ and considering the corresponding structure constants one can form an $E$-division algebra $A'$ which is $E$-isomorphic to $A$ such that $F'$ is a maximal subfield of $A'$.

We shall show that $A' \tensor_E E'$ is an $E'$-division algebra, but first let us show that this implies that $F'E'/E'$ is an adequate $G$-extension (and hence $G$ is $E'$-admissible). Indeed, if $A' \tensor_E E'$ is a division algebra, then $F'\otimes_E E'$ is a field. It follows that $F'\otimes_E E'\cong F'E'$, since $F'\otimes_EE'$ is $G$-Galois over $E'$ [Sal99, Theorem 6.3]. Thus, $[F'E':E']=[F':E]$ and $F'\cap E'=E$. Since $F'E'$ splits $A' \tensor_E E'$ and as $\ind(A' \tensor_E E')=[F'E':E']$, $F'E'$ is a maximal subfield of $A' \tensor_E E'$ and hence an adequate $G$-extension.

In order to show that $A' \tensor_E E'$ is an $E'$-division algebra, we first note that for each $i \in I$, $P_i = F_iQ_i = F'Q_i$, by [HaV96, Lemma 3.6(b)]. Thus, we have the following diagram:

$$ 
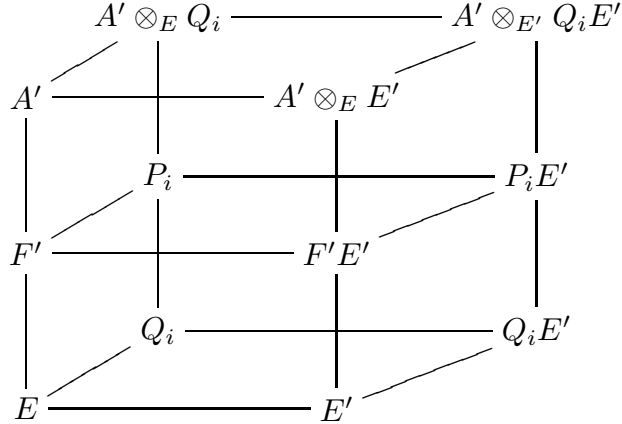
\xymatrix@=10pt{
& A' \otimes_E Q_i \ar@{-}[dl] \ar@{-}[rr] \ar@{-}[dd] && A' \otimes_{E'} Q_iE' \ar@{-}[dl] \ar@{-}[dd] \\
A' \ar@{-}[rr] \ar@{-}[dd] && A' \otimes_E E' \ar@{-}[dd]\\
& P_i \ar@{-}[dl] \ar@{-}[dd]\ar@{-}[rr] && P_iE' \ar@{-}[dl] \ar@{-}[dd] \\
F' \ar@{-}[rr] \ar@{-}[dd] && F'E'  \ar@{-}[dd] \\
& Q_i \ar@{-}[dl] \ar@{-}[rr] && Q_iE' \ar@{-}[dl] \\
E \ar@{-}[rr] && E'
}
$$

As mentioned above, $A'\otimes_E Q_i$ is Brauer equivalent to $D_i=D'_i\otimes_E Q_i$. Thus $A'\otimes_E Q_iE'$ is Brauer equivalent to $D'_i\otimes_{Q_i} Q_iE'$ which by choice of $D_i'$ is a division algebra. Then, $${n\over n_i}=\ind(D'_i)=\ind(D'_i\otimes_{Q_i} Q_iE')|\ind(A'\otimes_E Q_iE')|\ind(A'\otimes_E E')$$ for all $i\in I$. It follows that $n=\lcm_{i\in I}({n\over n_i}) | \ind(A'\otimes_E E')$. Hence $n=\ind(A'\otimes E')$, which shows that $A'\otimes_E E'$ is a division algebra. \qed

As a corollary, we get our main theorem.

\proclaim{Theorem \cp2.8} Let $R$ be a complete local domain of dimension $2$, with a separably closed residue field. Let $E = \Quot(R)$ and let $G$ be a finite group of order not divisible by $\chr(E)$. Then $G$ is $E$-admissible if and only if all the Sylow subgroups of $G$ are abelian of rank at most $2$.
\endproclaim \demo{Proof} By Proposition \cp2.7, if the Sylow subgroups of $G$ are abelian of rank at most $2$ then $G$ is $E$-admissible. For the converse, assume $G$ is $E$-admissible. For a prime $v$ of $E$, let $\ram_v$ denote the ramification map $\ram_v \colon \Br(E)\rightarrow \H^1(G_{E_v},\mQ/\Z) $ (see [Sal99]). Following [HHK10], we say that an $\alpha \in \Br(E)$ is determined by ramification with respect to a set of primes $\Omega$ if there is a prime $v\in \Omega$ for which $\exp(\alpha)=\exp(\ram_v(\alpha))$. Let $D$ be an $E$-division algebra with maximal subfield $L$ that has Galois group $G=\Gal(L/E)$. By [HHK10, Theorem 3.3], if $D$ satisfies:

\noindent (1) the order of $D$ is prime to $p$
 and $\ind(D)=\exp(D)$,

\noindent(2) $D$ is determined by ramification with respect to some set of discrete valuations,

\noindent then $G$ has Sylow subgroups that are abelian of rank at most $2$. Condition (1) is satisfied for any $\alpha$ of order prime to $p$ by [COP02, Theorem 2.1], while Condition (2) is satisfied by [COP02, Corollary \cp1.9 (c)] with respect to the set of codimension $1$ primes of $R$. \qed

\rem{Remark 2.9} Let $E$ be as above. By [COP02, Theorem 2.1], any Brauer class $\alpha\in \Br(E)$ of order prime to $\chr(E)$ has $\ind(\alpha)=\exp(\alpha).$ Thus by [Sch68, Proposition 2.2], a subfield of an $E$-division algebra is also a maximal subfield of some $E$-division algebra. \endrem

\references {HaJ98a}




\ref [BHH10] L.~Bary-Soroker, D.~Haran and D.~Harbater, \sl
Permanence criteria for semi-free profinite groups, \rm Mathematische Annalen, to appear.

\ref [COP02] J.L.~Colliot-Th\'el\`ene, M.~Ojanguren and R.~Parimala, \sl Quadratic forms over function fields
of two-dimensional Henselian rings and Brauer groups of related schemes, \rm Proceedings of the International Colloquium on Algebra, Arithmetic and Geometry, Tata Inst. Fund. Res. Stud. Math., vol. 16, pp. 185--217, Narosa Publ. Co., 2002.

\ref [DI71] F.~DeMeyer and E.~Ingraham, \sl Separable algebras over commutative rings, \rm Springer-Verlag, Berlin, 1971.

\ref [Fei04] W.~Feit, \sl PSL2(11) is admissible for all number fields, \rm Algebra, arithmetic
and geometry with applications (West Lafayette, IN, 2000), Springer,
Berlin, 2004, pp. 295–-299.

\ref [FrJ05] M.~Fried and M.~Jarden, \sl Field Arithmetic, 2nd
edition, revised and enlarged by M. Jarden, \rm Ergebnisse der
Mathematik III {\bf 11}, Sprin\-ger 2005.

\ref [FrP04]
J.~Fresnel and M.~v.d.~Put,
\sl Rigid Analytic Geometry and its Applications,
\rm Progress in Mathematics {\bf 218}, Birkh\"auser Boston, 2004.


\ref [FS95a] B.~Fein and M.~Schacher, \sl \bbQ(t) and \bbQ((t))-admisisbility of groups of odd order, \rm Proc. Amer. Math. Soc. {\bf 123} (1995).

\ref [FS95b] B.~Fein and M.~Schacher, \sl Crossed products over algebraic function fields, \rm Journal of Algebra {\bf 171} (1995), 531--540.

\ref [FSS92] B.~Fein, D.~Saltman and M.~Schacher, \sl Brauer-Hilbertian fields, \rm Trans. Amer. Math. Soc.
{\bf 334} (1992), 915--928.

\ref [Ha87] D.~Harbater,
\sl Galois coverings of the arithmetic line,
\rm in Number Theory --- New York 1984--85, ed. by D. V. and
G. V. Chudnovsky, Lecture Notes in Mathematics {\bf 1240},
Springer, Berlin, 1987, pp. 165--195.

\ref [HaH09] D.~Harbater and J.~Hartmann, \sl Patching over fields, \rm Israel Journal of Mathematics, to appear.

\ref [HaJ98] D.~Haran and M.~Jarden, \sl Regular split embedding
problems over complete valued fields, \rm Forum Mathematicum {\bf
10} (1998), 329--351.

\ref [HaS05] D.~Harbater and K.~Stevenson, \sl Local Galois theory
in dimension two, \rm Advances in Mathematics {\bf 198} (2005),
623--653.

\ref [HaV96] D.~Haran and H.~V\"olklein, \sl Galois groups over
complete valued fields, \rm Israel Journal of Mathematics {\bf 93}
(1996), 9--27.

\ref [HHK10] D.~Harbater, J.~Hartmann and D.~Krashen, \sl Patching subfields of division algebras, \rm to appear.
%




%
%
%

\ref [Jac96] N.~Jacobson, \sl Finite-dimensional division algs. over fields, \rm
Springer-Verlag, Berlin, 1996.

\ref [Lid96] S.~Liedahl, \sl K-admissibility of wreath products of cyclic p-groups, \rm J. Number Theory {bf 60} (1996), no. 2, 211--232

\ref [Mat86] H.~Matsumura, \sl Commutative ring theory, \rm
Cambridge University Press, 1986.

\ref [Nef10] D.~Neftin, \sl Admissibility and realizability, \rm submitted. arXiv:0904.3772v1.

%
%


\ref [Pa08] E.~Paran, \sl Algebraic patching over complete
domains, \rm Israel Journal of Mathematics {\bf 166} (2008),
185--219.

\ref [Pa10] E.~Paran, \sl Galois theory over complete local domains
domains, \rm to appear.

\ref [Po96] F.~Pop, \sl Embedding problems over large fields, \rm
Annals of Mathematics {\bf 144} (1996) 1--34.

\ref [Po10] F.~Pop, \sl Henselian implies large, \rm
to appear.

\ref [Rei75] I.~Reiner, \sl Maximal orders, \rm Academic Press, New York, 1975.

\ref [Sal99] D.J.~Saltman, \sl Lectures on division algebras,
\rm CBMS Regional Conference series no. 94. American Mathematical society, Providence RI, 1999.

\ref [Sch68] M.~Schacher, \sl Subfields of division rings, \rm Journal of Algebra {\bf 9} (1968), 451-–477.

\ref [Son83] J.~Sonn, \sl Q-admissibility of solvable groups, \rm Journal of Algebra {\bf 84} (1983), no. 2,
411–-419.

\ref [SS92] M.~Schacher and J.~Sonn, \sl K-admissibility of $A_6$ and $A_7$, \rm Journal of Algebra
{\bf 145} (1992), no. 2, 333–-338.

\ref [Ste82] L. Stern, Ph.D Thesis. Technion 1982.
%
%
%
%

\bye